\documentclass{article}
\newtheorem{thm}{Theorem}[section]
\newtheorem{prop}[thm]{Proposition}
\newtheorem{cor}[thm]{Corollary}
\newtheorem{lem}[thm]{Lemma}
\newtheorem{de}[thm]{Definition}
\newtheorem{rem}[thm]{Remark}

\newtheorem{ob}[thm]{Observation}
\newenvironment{proof}{\textbf{Proof.}}{\hfill\qed}
\def\sp{\vskip 3 true mm \noindent}
\def\wt{\widetilde}
\def\p{\par\noindent}

\def \wtA{\widetilde{A}}

\def \wtE{\widetilde{E}}
\def \wtF{\widetilde{F}}

\def \wtK{\widetilde{K}}

\def \wta{\widetilde{a}}
\def \wtb{\widetilde{b}}
\def \wte{\widetilde{e}}

\def \wtx{\widetilde{x}}
\def \wty{\widetilde{y}}
\def \wtz{\widetilde{z}}

\def \ot{\p \hbox to 8 cm{\hfill}}
\def \otb{\p \hbox to 2cm{\hfill}}

\def \r{\rightarrow}

\def\mapright{\mathop{\longrightarrow}\limits}
\def\mapleft{\mathop{\longleftarrow}\limits}
\def\ol{\overline}

\def\Xul{\underline{X}}
\def\aol{\overline{a}}
\def\bol{\overline{b}}

\def\pol{\overline{p}}

\def\Dol{\overline{D}}
\def\Fol{\overline{F}}

\def\Pol{\overline{P}}

\newcommand{\K}{\mathbb{K}}

\newcommand{\Z}{\mathbb{Z}}

\def \5{[\hskip 5mm]}
\def\sqr#1#2{{\vcenter{\vbox{\hrule height.#2pt
        \hbox{\vrule width.#2pt height#1pt \kern#1pt
	   \vrule width.#2pt}
        \hrule height.#2pt}}}}
\def\dqed{\mathchoice\sqr65\sqr65\sqr{5.1}4\sqr{4.5}4}
\def\qed{$\dqed$}
\font\msbmten=msbm10 		
\font\msbmseven=msbm7		
\font\msbmfive=msbm5		

\newfam\amssym
\textfont\amssym=\msbmten	
\scriptfont\amssym=\msbmseven
\scriptscriptfont\amssym=\msbmfive

\begingroup
\multiply\count0 by "100\advance \count0 by '157
\global\mathchardef\rtimes=\count0
\endgroup

\newcounter{enumrom}
\newenvironment{enumrom}{\begin{list}{\roman{enumrom})}{\usecounter{enumrom}
\setlength{\rightmargin}{\leftmargin}}}{\end{list}}
\usepackage[all,v2]{xy}
\LaTeXdiagrams
\usepackage{amssymb}
\title{Good Reduction of Good Filtrations at Places}
\author{Toukaiddine Petit\footnote{{\tt Author supported by the Scientific Programme NOG of the European Science Foundation}. 
}\\e-mail: toukaiddine.petit@ua.ac.be\\ Fred Van Oystaeyen\footnote{{\tt Acknowledging the EC project Liegrits MCRTN 505078}. 
}\\e-mail: Fred.vanoystaeyen@ua.ac.be
\\Department Wiskunde en Informatica, Universiteit Antwerpen\\
B-2020 (Belgium)}   
\date{}

\begin{document}
\maketitle
\subsubsection*{Abstract}
We consider filtered or graded algebras $A$ over a field $K$. Assume that there is
a discrete valuation $O_v$ of $K$ with $m_v$ its maximal ideal and $k_v:=O_v/m_v$ its residue field. Let $\Lambda$ be $O_v$-order such that $\Lambda K=A$ and $\overline{\Lambda}:=k_v\otimes_{O_v}\Lambda$ the $\Lambda$-reduction of $A$ at the place $K\leadsto k_v$. As in many examples of quantized algebras $A$ comes with a specific filtration that reduces well with respect to the valuation filtration defined by $\Lambda$ on $A$ and the reduction relates to the part of degree zero in the associated graded algebra. Hence several lifting properties fellow from valuation like theory, also for modules with good filtrations.
\sp
\noindent \textbf{Classification AMS 2000}: 16W35, 16W70, 16W60, 06B23, 06B25.
\sp
\noindent \textbf{Key words}: filtered algebras, valuations, reductions, quantum groups, generalized Weyl algebras.
\section*{Introduction}
One possible arithmetical aspect in the noncommutative geometry of associative
algebras may be found in the construction of a noncommutative divisor theory
based on noncommutative valuations, e.g. \cite{VO}.  Reduction of algebras at
such valuations have already been investigated in (\cite{HVO1}, \cite{MVO}).  Typical
algebras considered there are among others : rings of differential operators,
certain quantum groups, quantized algebras and regular algebras in the sense
of projective noncommutative algebraic geometry.  these algebras have a
natural gradation or filtration defined in terms of some finite dimensional
vector spaces, e.g. the part of degree one is finite dimensional.  In this
note we study the reduction of the filtered or graded structures over a given
valuation in the base field, $K$ say.  Its properties relate to
certain lattices in the characteristic vector spaces hinted at above.  For some
filtration $FA$ on a $K$-algebra $A$ the unramifiedness property of a
reduction relates to the induction of good filtrations (cf. \cite{Bj}) in every
$F_nA$.  Perhaps the main result in this context is the establishing of a
lifting property for unramified reductions from the associated graded ring
$G_F(A)$ to the filtered ring $A$.
\sp
Several interesting classes of algebras may be studied via reduction
techniques.  The color Lie algebras and their enveloping algebras will be
separately treated in forthcoming work.  An important class of examples
consists of generalized Weyl algebras (cf. \cite{B1}) or generalized crossed products
(cf. \cite{BVO}); this class contains popular algebras like~: quantum deformed Weyl
algebras, the quantum plane, quantum $\mathcal{U}_q(sl_2)$ of $sl_2$, the quantum
Heisenberg algebra (cf. \cite{Ma}), Witten's first and Woronowicz's deformation,
the quantum group $O_{q^2}$ of $so_3$ (cf. \cite{Sm}) etc... .  For algebras in the
foregoing class the extension of valuations on the base field to noncommutative
valuations on their fields of fractions has been studied (cf. \cite{MVO}) and
several lifting results for regularity conditions as well as dimension
calculations follow from the reduction properties.  
\sp
As a general reference for detail on filtered rings and modules we refer to
\cite{HVO1}, full detail on graded ring theory may be found in (\cite{NVO1},\cite{NVO2}).

\section{Preliminaries on Reductions, Filtrations and Gradations}

Throughout $A$ is an associative algebra over a commutative field $K$.  A
$\Bbb{Z}$-filtration $FA$ is given by an ascending family $\{F_nA, n\in
\Bbb{Z}\}$ of additive subgroups such that $F_nAF_mA\subset F_{n+m}A$ for all
$n,m\in\Bbb{Z},1\in F_0A$, and we always assume the filtration to be
exhaustive, i.e. $A=\cup_{n\in\Bbb{Z}}F_nA$, and separated, i.e.
$0=\cap_{n\in\Bbb{Z}}F_nA$. We say that $A$ is a filtered $K$-algebra of
$K\subset F_0A$, consequently all $F_nA$ are $K$-vector spaces.  Following
conventions and notation of (\cite{HVO1}, \cite{NVO1}), we write $G_F(A)$ for the
associated graded ring, or $K$-algebra, with respect to $FA$ and we let
$\widetilde{A}$ be the Rees ring or blow-up ring, respectively $K$-algebra.
We write~: $G_F(A)=\oplus_{n\in\Bbb{Z}}G_F(A)_n$ with $G_F(A)_n=F_nA/F_{n-1}A$
for all $n\in\Bbb{Z}$, $\widetilde{A}=\oplus_{n\in\Bbb{Z}}\widetilde{A}_n$ with
$\widetilde{A}_n=F_nA$ for $n\in\Bbb{Z}$.  It is practical to identify
$\wt{A}$ with the graded subring $\sum_{n\in\Bbb{Z}}F_nAT^n$ in $A[T,T^{-1}]$
where $T$ is a central variable of degree one.  Recall that the so-called
principle symbol map $\sigma_F:A\r G_F(A)$ is defined by mapping an $a\in A$
such that $a\in F_nA-F_{n-1}A$ to $a\ \mathrm{mod}\ F_{n-1}A$ in $G_F(A)_n$; observe
that $\sigma_F$ is neither addictive nor multiplicative in general. If no
ambiguity can arise the subscript may be dropped in notation introduced above.
\sp
Zariskian filtrations on noncommutative rings have been characterized in
several ways (cf. \cite{HVO1}) but in any case these filtrations have the property
that $A,G_F(A)$ and $\wt{A}$ are (twosided) Noetherian rings.  
\sp
A filtration on a $K$-algebra $A$ is said to be a 
{\bf finite filtration} if ${\rm dim}_KF_nA$ is finite for all $n\in\Bbb{Z}$.
Similarly, a graded algebra $R=\oplus_nR_n$ is said to be {\bf finitely
graded} if ${\rm dim}_KR_n$ is finite for all $n\in \Bbb{Z}$.  Obviously, if
$FA$ is finite then $G_F(A)$ and $\wt{A}$ are both finitely graded; if
$\wt{A}$ is finitely graded then $FA$ is finite and $G_F(A)$ is finitely
graded.  If $G_F(A)$ is finitely graded then $FA$ is finite if and only if at
least one $F_mA$ is finite dimensional over $K$.  Typical graded algebras
appearing in noncommutative projective geometry e.g. regular algebras as
studied in \cite{Art}, \cite{ATV}, are graded $K$-algebras of type $R=K\oplus R_1\oplus
\ldots$, generated by $R_1$ over $K$ as a $K$-algebra and ${\rm dim}_KR_1$
being finite dimensional.
\sp
Let us recall some definitions and facts concerning valuations of skewfields,
the old book of O. Schilling is still a valid basic reference for the general
theory, cf. \cite{Sc}.  A subring $\Lambda$ in a skewfield $\Delta$ is said to be
{\bf a valuation ring} of $\Delta$ if for every $x\in\Delta-\{0\}$ either $x$
or $x^{-1}$ is in $\Lambda$ and moreover $\Lambda$ is invariant under inner
automorphisms of $\Delta$.  The unique maximal ideal $P$ of $\Lambda$ given by
$P=\{x\in \Lambda,x^{-1}\not\in \Lambda\}$ defines the residue field (!)
$\Lambda/P$ of $\Delta$; we often write $\Delta_v=\Lambda/P$ (sometimes
$\ol{\Delta}=\Lambda/P$).  A valuation ring $\Lambda$ of $\Delta$ is said to
be {\bf discrete} if $P$ is a principal ideal or equivalently $\Lambda$ is
Noetherian and the value group is $\Bbb{Z}$.  When $\Delta$ is a $K$-algebra
and $\Lambda$ is a valuation ring of $\Delta$ then $\Lambda\cap K$ is a
valuation ring of $K$; in case $K\subset \Lambda$ we say that $\Lambda$ is a
$K$-valuation ring.  We write $O_v\subset K$ for a valuation ring of $K$ and
denote its maximal ideal by $m_v$ and its residue field by $k_v=O_v/m_v$.
From a valuation ring $O_v\subset K$ we derive a valuation function $v:K^*\r
\Gamma$ for a suitable totally ordered abelian group; in the discrete case we
are looking at $\Gamma=\Bbb{Z}$. To a noncommutative valuation ring $\Lambda$
in $\Delta$ we may also associate a valuation function $\nu:\Delta^*\r \Gamma$
where now $\Gamma$ is again totally ordered but not necessarily
abelian.  In some cases the abelian property of $\Gamma$ is forced upon us,
e.g. noncommutative valuation of the skewfield of the first Weyl algebra are
necessarily having an abelian value group.  In the sequel, unless otherwise
stated, all {\bf valuation are supposed to be discrete} e.g. in particular we
only consider $\Bbb{Z}$-valuations.  If $\Lambda$ is a noncommutative discrete
valuation ring of $\Delta$ then we define a filtration $F^v\Delta$ on
$\Delta$, called the {\bf valuation filtrations}, by putting
$F^v_n\Delta=P^{-n}$.  If $\Lambda=O_v,\Delta=K$ then we write $f^vK$
for the valuation filtration of $K$.   Observe that~: ${\rm deg}\sigma_{F^v}
(\delta)=-v(\delta)$ for $\delta\in\Delta$, ${\rm deg}_{f^v}(x)=-v(x)$ for
$x\in K$.  In the situation $K\subset \Delta$ and for a given valuation ring
$\Lambda$ of $\Delta$ with valuation function $\nu$ the valuation ring
$\Lambda\cap K$ of $K$ is the {\bf induced valuation ring}, denoted by $O_v$.
Of course $P\cap K=m_v$ but it is possible that $P^e\cap K=m_v$ for $e>1$.
Since $\cap _{n\in\Bbb{N}}P^n=0$ it follows that there is a unique $e_{\nu}$
such that $\pi\in P^{e_{\nu}}$ but $\pi\not\in P^{e_{\nu}+1}$ where
$m_v=(\pi)\subset O_v$.  This $e_{\nu}\in\Bbb{N}$ is called the {\bf
ramification index of $\Lambda$ over $O_v$}.  We easily check that $P^m\cap
K=m^d_v$ where $d=\lceil{m\over e}\rceil$ is the smallest integer bigger
that or equal to ${m\over e}$.  This shows that $e_{\nu}$ is in fact the
ramification of the valuation filtration $F^v\Delta$ over $f^vK$, i.e.
$F^v_m\Delta\cap K=f^v_dK$ where $d=\lceil{m\over e}\rceil$ as above.
Whereas the $m_v$-adic filtration of $\Lambda$ obviously induces $f^vK$
or in fact the negative part of it viewed as a filtration on $\Lambda$, such
statement is false for $F^v$ as noted before.  In general a filtration $FR$ of
ring $R$ is said to be {\bf scaled with step $d$} if $F_0R=R_1R=\ldots\subset
F_dR=F_{d+1}R=\ldots \subset F_{2d}R=F_{2d+1}R=\ldots$ and similar on the
negative side~: $\ldots\subset F_{-d}R\subset
F_{-d+1}R=\ldots=F_{-1}R\subset_0R$, i.e. $FR$ is obtainable as a
$d\Bbb{Z}$-filtration viewed as a $\Bbb{Z}$-filtration in the way explained
above.
\begin{lem}\label{l1.1}
With notation as before, $F^v\Delta$ induces in $K$ the scaled filtration with
step $e$ associated to $f^vK$, where $e$ in the ramification index of
$\Lambda$ over $O_v$.
\end{lem}
\sp
In the foregoing it is obvious that $m_v$ is contained in the Jacobson radical
$J(\Lambda)$ of $\Lambda$; this is so because we assumed that $\Lambda$ is a
valuation ring extending $O_v$ e.g. $P\cap K=m_v$.
Observe however that for an arbitrary $O_v$-order $\Lambda$ in an infinite
dimensional $\Delta$ over $K$ we need not have $m_v$ in the Jacobson
radical of $\Lambda$, the latter may even be zero (verify for the Weyl algebra
defined over $\Bbb{Z}_p$ as an order in the Weyl field over $Q$)~!
\sp
Let us recall Proposition 3.1. from \cite{HVO1}.
\begin{prop}\label{p2.1}
Let $R$ be an Artinian ring with filtration $FR$, then the following
statements are equivalent~:
\begin{enumrom}
\item
$G_F(R)$ is a domain
\item
$R$ is a skewfield and every nonzero homogeneous element of $G_F(R)$ is
invertible i.e. $G_F(R)$ is a graded-skewfield.
\item
$R$ is a skewfield, $F_0R$ is a discrete valuation ring of $R$ with maximal
ideal $F_{-1}R$ and $F_{-ne}R=(F_{-1}R)^n$ for some $e\in\Bbb{N}$.  
\end{enumrom}
\end{prop}
\section{Reductions of Gradations and Filtrations}
Again we look either at (separated and exhaustive) filtered $K$-algebras $A$
with a subring $\Lambda$ such that $\Lambda\cap K=O_v$, or else at graded
$K$-algebras $\wtA$ with a subring $\wt{\Lambda}$ that is a graded subring now
such that $\wtA\cap K=O_v$.  In the sequel we shall only consider $O_v$-orders
$\Lambda$, resp. $\wt{\Lambda}$, such  that $K\Lambda=A$, resp.
$K\wt{\Lambda}=\wtA$.  So we have the induced filtration $F\Lambda$ given by
$F_n\Lambda=\Lambda \cap F_nA$, or the induced gradation
$\wt{\Lambda}_n=\Lambda\cap \wtA_n$.
\begin{ob}\label{o2.1}
With notation as before~:
\begin{enumrom}
\item
$m^a_v\Lambda\cap F_nA=m^a_{v}(\Lambda\cap F_nA)$, for all $n\in \Bbb{Z}$,
$a\in\Bbb{Z}$.
\item
$m^a_v\wt{\Lambda}\cap\wtA_n=m^a_v(\wt{\Lambda}\cap\wtA_n)$, for all
$n\in\Bbb{Z}, a\in\Bbb{Z}$.
\end{enumrom}
\end{ob}
\begin{proof}
Let us establish i), the proof of ii) is similar.
\begin{enumrom}
\item
The inclusion $m^a_v(\Lambda\cap F_nA)\subset m^a_v\Lambda\cap F_nA$ is
trivial.  Pick $z\in m^a_r\Lambda\cap F_nA$, i.e. $z=\pi^a\lambda$ for some
$\lambda\in\Lambda$.  Since $F_nA$ is a $K$-space $\lambda=\pi^{-a}Z\in
F_nA\cap \Lambda$, hence $z\in \pi^a(F_nA\cap \Lambda)$ and because $F_nA\cap
\Lambda$ is an $O_v$-module the latter equals $m^a_v(F_nA\cap
\Lambda)$.
\end{enumrom}
\end{proof}

In the situation as above we call $\Lambda/m_v\Lambda$ the {\bf
$\Lambda$-reduction of $A$ at the place
$K\leadsto k_v$}
(or at $v$); similarly $\wt{\Lambda}/m_v\wt{\Lambda}$ is the
$\Lambda$-reduction of $\wt{A}$ and it is clearly a graded $k_v$-algebra.
Faithful to the notation of the residue field we write
$A_v=\Lambda/m_v\Lambda$, $\wtA_v=\wt{\Lambda}/m_v\ol{\Lambda}$ and we write
$\pi:\Lambda\r A_v,\wt{\pi}:\wt{\Lambda}\r \wt{A}_v$ for the corresponding
canonical ring epimorphisms.  From the observation i) it is clear that $FA$
defines a filtration $FA_v$ (in fact i expresses a compatibility relation
between $FA$ and $F^vA$ defined by $F^v_uA=m^{-n}_v\Lambda$ i.e. the $m_v$-adic
$\Bbb{Z}$-filtration of $A$ constructed from $\Lambda$) given by
$F_nA_v=F_n\Lambda/m_vF_n\Lambda$.  Moreover we have a graded subring
$G_F(\Lambda) \subset G_F(A)$ such that $G_F(\Lambda)\cap K=O_v$
and $KG_F(\Lambda)=G_F(A)$, so this defines an exhaustive  (graded)
filtration $f^vG_F(A)$ by $f^v_nG_F(A)=m^{-n}_vG_F(\Lambda)$
On the other hand $G_{F^v}(A)=\Lambda\otimes_{o_v}k_v[t,t^{-1}]\cong
A_v[t,t^{-1}]$ has a filtration induced by $FA$, via $F\Lambda$,
let us denote it by
$\Fol G_{F^v}(A)$, then $\Fol_nG_{F^v}(A)\cong (F_n\Lambda/m_rF_n\Lambda)
[t,t^{-1}]$ (we may view this as an identification if we identify $G_{F^v}(A)$
and $A_v[t,t^{-1}]$.  The compatibility between $FA$ and $F^vA$ actually
establishes : $G_{f^v}(G_F(A))=G_{\Fol}(G_{F^v}(A))$.  In general we do
not know that the filtration $F^vA$ associated to $\Lambda$ is separated, but
when suitable finiteness conditions hold all filtrations constructed  before
will be separated.  let is first mention a different easy but sometimes
interesting good case.
\begin{lem} \label{l2.2}If $\Lambda$ and $A$ have no nonzero ideal in common then
$F^vA$ is separated.  A similar statement holds with respect to
$\wt{\Lambda}$ and $\wt{A}$ in the graded case.
\end{lem}
\begin{proof}
Put $E=\cap\{n\in \Bbb{N},m^n_v\Lambda\}\subset \Lambda$.  If $x\in E$ then
$\pi^{-n}x\in E$ for every $n\in\Bbb{N}$, thus $Kx\subset E$ and also
$Ax\subset E$, similar for $AxA\subset E$.  This leads to a contradiction were
$x\neq 0$.
\end{proof}
\begin{cor}\label{c2.3}
If $A$ is a simple $K$-algebra then $F^vA$ is always separated, for every
$O_v$-order $\Lambda$.
\end{cor}
\begin{prop} \label{p2.4} With notation as before, if $f^vG_F(A)$ is a
separated filtration then $F^vA$ is separated.
\end{prop}
\begin{proof}
In view of Observation \ref{o2.1} we have that
$f^v_nG_F(A)_d=m^{-n}_vG_F(\Lambda)_d$, for all $d\in\Bbb{Z}$, and
$m^{-n}_vG_F(\Lambda)_d=G_F(A)_d\cap m^{-n}_vG_F(\Lambda),
G_F(m^{-n}_v\Lambda)=m^{-n}_vG_F(\Lambda)$.  Hence $G_F(E)\subset
\cap_{n\in\Bbb{N}}m^n_vG_F(\Lambda)=0$, the latter following from the assumed
separateness of $f^vG_F(A)$.  Thus $E\subset\cap_{n\in\Bbb{N}}F_nA$ but as
$FA$ is separated (that was a standing assumption throughout) it follows that
$E=0$, hence $F^vA$ is separated too.
\end{proof}
\begin{de}  We say that $\Lambda$ is {\bf $FA$-finite} if for
all $d\in\Bbb{Z}, \Lambda_d=\Lambda\cap F_dA$ is a finitely generated $O_v$-module.  In
the graded situation $\wt{\Lambda}\subset \wt{A}$ we say that {\bf
$\wt{\Lambda}$ is $\wt{A}$-finite} if $\wt{A}_d\cap
\wt{\Lambda}=\wt{\Lambda}_d$, for all $d\in\Bbb{Z}$, is a finitely generated
$O_v$-module.  For a finite dimensional vetorspace $V$ over $K$, an
$O_v$-module $M$ contained in $V$ is said to be an {\bf $O_v$-lattice} of $V$
if ${\rm rank}_{O_v}M={\rm dim}_KV$.  Any $O_v$-lattice $M$ of $V$ defines an
{\bf unramified reduction} $V_v=M/m_vM=k_v\otimes_{O_v}M$ with ${\rm
dim}_{k_v}V_v={\rm dim}_KV$.
\end{de}
\begin{thm}\label{t2.6}
With notation and conventions as before~:

\begin{enumerate}
	\item If $G_F(\Lambda)$ is $G_F(A)$-finite then $F^vA$ is separated.
	\item 
\begin{enumerate}
	\item If $\Lambda$ is $FA$-finite then $f^vG_F(A)$ and $F^vA$ are both separated\\
filtrations.  The restriction of $\Fol G_{F_v}(A)$ to $A_v=G_{F^v}(A)_0$
denoted by $\Fol A_v$, is given by $\Fol_nA_v=F_n\Lambda/m_vF_n\Lambda$
which is an unramified reduction of $F_nA$.  Moreover
$G_{F^v}(A)=A_v[t,t^{-1}]$ and it has the residual filtration given by,
$\Fol_nA_n[t,t^{-1}]$ $n\in \Z$
	\item The filtration $F^vA$ induces a good filtration in $F_dA$ for every
$d\in\Bbb{Z}$.
\end{enumerate}
\end{enumerate}
\end{thm}
\begin{proof}
1. If $G_F(\Lambda)$ is $G_F(A)$-finite then $G_F(A)$ is a finite graded
$K$-algebra; since every $G_F(\Lambda)_d, d\in \Bbb{Z}$, is a finitely
generated and torsion free $O_v$-module it is free of rank $n_d$.  Now $E_{\rm
gr}=\cap_{n\in\Bbb{N}}m^n_vG_F(\Lambda)$ is a common graded ideal of
$G_F(\Lambda)$ and $G_F(A)$ with $R_{{\rm gr},d}\subset G_F(\Lambda)_d$, the
latter free of finite rank $n_d$ over $O_v$.  As observed for $E$ earlier,
also for $E_{\rm gr}$ we do have that $\pi E_{\rm gr}=E_{\rm gr}$ and since
$\pi\in G_F(\Lambda)_0$ we also have $\pi E_{{\rm gr},d}=E_{{\rm gr},d}$ for
every $d\in\Bbb{Z}$.  Since now we are dealing with finitely generated
$O_v$-modules Nakayama's lemma yields that $E_{{\rm gr},d}=0$ for all $d\in
\Bbb{Z}$, hence $E_{\rm gr}=0$ or $f^vG_F(A)$ is separated.  Foregoing
proposition \ref{p2.4} then yields that $F^vA$ is separated.\\
2.a. If $\Lambda$ is $FA$-finite then $G_F(\Lambda)$ is $G_F(A)$-finite, hence
$f^vG_F(A)$ and $F^vA$ are both separated.  In view of the finiteness
assumption $F_{\alpha}\Lambda$, for every $d\in\Bbb{Z}$, is an $O_v$-lattice
hence a free $O_v$-module of rank $n_d$ say; then
$F_d\Lambda/m_vF_{\alpha}\Lambda$ is a $k_v$-vector space of dimension $n_d$,
thus $\Fol_nA_v$ is indeed an unramified reduction of $F_nA$.  The remaining
claims are just reformulations of earlier observations.\\
2.b. Recall that for a filtered modules $M$, with filtration $FM$, over the
filtered ring $A$ we say that $FM$ is a good filtration if there is a finite
set $m_1,\ldots,m_s$ in $M$ such that for every $n\in\Bbb{Z}$ we have that~:
$F_nM=\sum^s_{i=1}F_{n-d_i}A.m_i$, for some fixed $d_1,\ldots,d_s$ in
$\Bbb{Z}$.  Now from $\Lambda m^p_v\cap F_nA=m^p_vF_n\Lambda$ for all $n,p$,
it is clear that the filtration induced in $F_nA$ is good (viewed as a
filtered $f^vK$-module). Indeed it suffices to pick an $O_v$-basis for
the free $O_v$-module $F_n\Lambda$ for the $m_i$ and take each $d_i$ to be
zero, then the only way to express an element of $m_vF_n\Lambda$ in the
selected basis is by taking coefficients from $m_v$.

\end{proof}

Now we look at a graded $O_v$-order $\wt{\Lambda}$ in $\wt{A}$ as before and
we assume that $\wt{\Lambda}$ contains a central regular homogeneous element
of degree one, $T$ say.  Put $A=\wt{A}/\wt{A}(T-1)$,
$\Lambda=\wt{\Lambda}/\wt{\Lambda}(T-1)$; then $A$ has a filtration $FA$ given
by $F_nA=\wtA_n/\wtA(T-1)\cap\wtA n$, and $\Lambda$ has a filtration
$F\Lambda$ given by $F_n\Lambda=\wt{\Lambda}_n/\wt{\Lambda}(T-1)\cap
\wt{\Lambda}_n$, for all $n\in\Bbb{Z}$.

\begin{lem}  With notation as before we obtain~:
\begin{enumrom}
	\item $\wtA(T-1)\cap \wt{\Lambda}=\wt{\Lambda}(T-1)$
	\item $F_n\Lambda=F_nA\cap\Lambda=\wt{\Lambda}_n/(T-1)\wt{\Lambda}\cap\wtA_n$.
\end{enumrom}
\end{lem}
\begin{proof}
i) Obviously $\wt{\Lambda}(T-1)\subset\wtA(T-1)\cap\wt{\Lambda}$.  For the
converse look at $\wt{z}(T-1)\subset\wt{\Lambda}$ with $\wt{z}\in\wt{A}$.  If
$\wtz\not\in\wt{\Lambda}$ then there is a $d\in\Bbb{Z}$ minimal such that
$\wt{z}_d\not\in\wt{\Lambda}_d$ while on the other hand~:
$\wtz_{d-1}T-\wtz_d\in\wt{\Lambda}_d$.  In case $\wtz_{n-1}\neq 0$, then the
foregoing entails that $\wt{z}_{d-1}T\not\in\wt{\Lambda}_d$ and thus
$\wt{z}_{d-1}\not\in\wt{\Lambda}_{d-1}$ because $T\in \wt{\Lambda}_1$ but that
contradicts minimality of $d$.  Hence $\wt{z}_{d-1}=0$, if $\wtz_{d-2}\neq 0$
then $\wtz_{d-2}T^2-\wtz_d\in\wt{\Lambda}$, thus as in the first part
$\wtz_{d-2}T^z\not\in\wt{\Lambda}_d$ and certainly $\wtz_{z-2}\not\in
\wt{\Lambda}_{d-2}$ again contradicting minimality of $d$.  So we are in the
situation where $\wtz_d$ is the homogeneous part of lowest degree in the
decomposition of $\wtz$.  They from $\wtz(T-1)\in\wt{\Lambda}$ we obtain that
the homogeneous part of lowest degree in the decomposition of $\wtz(T-1)$, and
that is exactly $-\wtz_d$, must be in $\wt{\Lambda}$ and that leads to a
contradiction.  Consequently $\wtz(T-1)\in\wt{\Lambda}$ leads to
$\wtz\in\wt{\Lambda}$ and the claim i. follows.\\
ii) From i) it is clear that $\Lambda\subset A$ and $F_n\Lambda\subset F_nA$ for
all $n\in\Bbb{Z}$.  If $a_n\in F_nA\cap \Lambda$ then there exists a
$\wt{\lambda}\in\wt{\Lambda}$ such that $\wt{\lambda}$ mod $(\wt{\Lambda}\cap
\wtA(T-1))=a_n$ but also there is an $\wt{a}_n\in\wtA$ such that $\wta_n{\rm
mod}(\wtA_n\cap \wtA(T-1))=a_n$.  Thus $\wta_n+\wtb(T-1)=\wt{\lambda}$ for
some $\wtb\in\wtA$, yields~: $(*)\wt{\lambda}_n=\wta_n+\wtb_{n-1}T-\wtb_n$.
Also $\wt{\lambda}_{n-1}=\wtb_{n-2}T-\wtb_{n-2}T-\wtb_{n-1}$ or
$\wtb_{n-1}T=\wtb_{n-2}T^2-\wt{\Lambda}_{n-1}T$.  Substituting an (*) then
leads to~: $\wt{\lambda}_n+\wt{\lambda}_{n-1}T=\wta_n+\wtb_{n-2}T^2$.  If
$\wtb_{n-2}\neq 0$ then we look at $\wt{\lambda}_{n-2}=\wtb_{n-3}T-\bol_{n-2}$
and arrive at
$\wt{\lambda}_n+\wt{\lambda}_{n-1}T+\wt{\lambda}_{n-2}T^2
=\wta_n+\wtb_{n-3}T^3$. We repeat this procedure until we obtain
$\wta_n=\wt{\lambda}_n+\wt{\lambda}_{n-1}T+\ldots \wt{\lambda}_{n-d}T^d$ and
thus $\wt{a}_n\in \wt{\Lambda}_n$.  From i. again it is clear that $a_n\in
F_n\Lambda$ follows, and the second equality of ii. also follows.

\end{proof}

Returning to the situation of $\Lambda\in A$ with filtration $FA$ inducing
$F\Lambda$, then the Rees ring (blow-up ring) of $A$ with respect to $FA$,
resp. $\Lambda$ with respect to $F\Lambda$, will be denoted by $\wtA$, resp.
$\wt{\Lambda}$.  Applying the foregoing lemma to these graded rings we recover
the filtered situation from the Rees ring situation.  Now the general theory
of filtered rings yields that $\wtA/T\wtA\cong G_F(A)$,
$\wt{\Lambda}/T(\wt{\Lambda})=G_F(\Lambda)$, with the graduation of $\wt{A}$,
resp, $\wt{\Lambda}$, defining the gradation of $G_F(A)$, resp. $G_F(\Lambda)$.
Let us write $\wtF^v$ for the graded filtration of $\wtA$ defined by
$\wtF^n_n\wtA=m^{-n}_v\wt{\Lambda}$.  The filtration $FA$ resp. $F\Lambda$,
corresponds to the $T$-adic filtration on $\wtA$, resp. $\wt{\Lambda}$ cf.
\cite{HVO1}.  We obtain the following extension of Proposition \ref{p2.4}:
\begin{lem}
If $\cap_n\wt{\Lambda}T^n=0$, then if $f^vG_F(A)$ is separated , recall
$f^v_nG_F(A)=m^{-n}_vG_F(\Lambda)$ for $n\in\Bbb{Z}$, then $\wtF^v\wtA$ is
separated too.
\end{lem}
\begin{proof}
Put $\wtE=\cap_{n\in\Bbb{Z}}\wtF^v_n\wtA=\cap_{n\in\Bbb{N}}m^n_v\wt{\Lambda}$.
Then (compare to Proposition \ref{p2.4}. , first part of proof) : $\wtE{\rm
mod}\wt{\Lambda}T\subset \cap_{n\in
\Bbb{N}}\pi^n\wt{\Lambda}/T\wt{\Lambda}=0$, i.e. $\wtE\subset T\wt{\Lambda}$.
Pick $\wte \in\wtE$, then we have $\wte=T\wt{\lambda}$ but also from
$\pi\wtE=\wtE$ it follows that $\wte=\pi\wtx$ for some $\wtx\in \wtE$ i.e.
$\wtx=T\wty$ for some $\wty\in\wtE$.  Consequently~: $\pi\wtx=\pi
T\wty=T\wt{\lambda}$ and since $T$ is regular in $\wt{\Lambda}, \pi\wty=\wt{\lambda}$ holds i.e. $\wt{\lambda}\in \wt{\Lambda}\pi$.  From
$e=T\pi\wt{\lambda}_1$ it follows that $T\wt{\lambda}_1=\wt{e}_1\in\wtE$ (note~:
$\pi^{-1}\wtE=\wtE$ !).  Repetition of the foregoing argument leads to
$\wt{\lambda}_1\in\wt{\Lambda}\pi$ etc... until we arrive at $\wte=\pi\wtx=T\wt{\mu}$ with $\wt{\mu}\in\wtE$, while $e=T^2\wt{\lambda}_1$ follows from
$\wt{\mu}=T\wt{\lambda}_1$ for some $\wt{\lambda}_1\in\wtE$, etc...  Finally
we obtain~:
$$
\wte=T\wt{\lambda}=T^2\wt{\lambda}_1=T^3\wt{\lambda}_2=\ldots=T^n\wt{\lambda}_n=\ldots \in \cap_{n\in\Bbb{N}}T^n\wt{\Lambda}=0$$
The proof is thus finished as $\wtE=0$ follows.
\end{proof}
\sp
The Rees ring of the valuation filtration $f^vK$ is $\wtK=O_v[\pi
t^{-1},\pi^{-1}t]$, where we now write $t$ for the regular homogeneous
element of degree one in $\wtK$.  The ring $\wtK$ is in fact a gr-valuation
ring in the field $K(t)$ of rational functions in $T$.  When calculating the
Rees ring of $A$ with respect to $F^vA$, $\wtA^{(v)}$ say, we may take 
$T=t\in \wt{A}^{(v)}_1$ and moreover $\wtA^{(v)}$ is a $\wtK$-algebra and it is 
strongly graded (recall that a graded ring $R$ is strongly graded if
$R_nR_{-n}=R_0$ for all $n\in\Bbb{Z}$, equivalently when $R_1R_{-1}=R_0$).
Note that the Rees ring of $A$ with respect to $FA$ need not have $t$ in
degree one, in fact one has to use another $T\in \wtA_1$ which relates to $t$
in some specific way reflecting the ramification of $FA$ over $f^vK$.
In particular $\wtA$ is not necessarily strongly graded (but it contains a
strongly $e\Bbb{Z}$-graded subring where $e$ is the ramification index of
$FA$ over $f^vK$).  We have a Rees version of Theorem \ref{t2.6}.
\begin{prop}\label{p2.9}
If $\Lambda$ is $FA$-finite then $\wt{\Lambda}$, respectively $A$, are defined
with respect to $F\Lambda$, respectively $FA$; the converse holds too.  Any of
the aforementioned properties entails that $G_F(\Lambda)$ is $G_F(A)$ finite.\\
If $G_F(\Lambda)$ is $G_F(A)$-finite then both $\wtF^{(v)}$ and $F^vA$ are
separated and so is $f^vG_F(A)$. \\
Conversely if $\wt{\Lambda}\subset \wtA$ are
given graded rings having a regular central homogeneous element of degree
one $T\in \wt{\Lambda}_1\subset \wtA_1$, then
$A=\wtA/(-T)\wtA,\Lambda=\wt{\Lambda}/(T)\wt{\Lambda}$ have filtrations $FA$,
resp. $F\Lambda$ such that $\wt{\Lambda}$, respectively $\wtA$, are indeed the Rees rings with
respect to those filtrations $F\Lambda$, respectively $FA$, and moreover
$G_F(\Lambda)=\wt{\Lambda}/T\wt{\Lambda},G_F(A)=\wtA/T\wtA$.  \\
The statements concerning (unramified) reductions as in Theorem \ref{t2.6} shift from filtered to
Rees level or back.
\end{prop}
\begin{proof}
All statements are consequence of earlier observations and results; let us just
point out that the property in Theorem \ref{t2.6}, 2.6.b., i,e, $F^vA$ inducing a good filtration
in $F_dA$, for every $d\in \Bbb{Z}$, viewed as a filtered $K$-module with
respect to $f^vK$, is just the finite generation property for the Rees
module of $F_dA$ with respect to $\wtK$ which is exactly
$\wtK\otimes_{O_v}F_d\Lambda=\wtF_d\wtA^{(v)}$.
\end{proof}
\sp
Let us finish this section by mentioning some further remarks about strong
filtrations, relating to valuations.  In general a filtration on a ring $A$ is
said to be a strong filtration if $F_nAF_mA=F_{n+m}A$ for all $n.m\in\Bbb{Z}$,
equivalently if $G_F(A)$ is a strongly graded ring, i.e.
$G_F(A)_nG_F(A)_m=G_F(A)_{n+m}$. By definition $F^vA$ is a strong filtration.
\begin{lem}\label{l2.10}
Let $A$ and $\Lambda$ be as before and consider an Ore set (left, right, left
and right) $S$ in $A$, then there is an Ore set $S_{\Lambda}$ in $\Lambda$
(left, right, let and right) such that $S^{-1}_{\Lambda}A=S^{-1}A$.
\end{lem}
\begin{proof}
Take $S_{\Lambda}=K^*S\cap \Lambda,K^*=K-\{0\}$.  It is straightforward to
check the (left, right, left and right) Ore conditions in $\Lambda$ for
$S_{\Lambda}$ and (obviously) $S^{-1}_{\Lambda}A=S^{-1}A$.
\end{proof}
\begin{cor}
The localized filtration derived from $F^nA$ is exactly the localized filtration
derived from the $m_v$-adic filtration of $\Lambda$; it is a strong filtration
denoted by $F^vS^{-1}A$.
\end{cor}
\sp
Let us say that $A$ is an {\bf order in an Artinian ring} if its set of
regular elements $S_0$ is a left and right Ore set such that $S^{-1}_0A$ is
Artinian (one sided left or right statements may be formulated similarly).
\begin{prop}\label{p2.12}
Let $A$ be an order in an Artinian ring $Q=S^{-1}_0A$ and $\Lambda\subset A$
as before.  If $A_v$ is a domain then $Q$ is a skewfield and $F^vA$ extends to
a strong filtration of $Q$ such that $F^v_0Q$ is a valuation ring of $Q$
extending $v$ from $K$ to $Q$.  Then in the situation where $\Lambda\subset A$
are graded, $A$ a graded $K$-algebra such that each $A_n$ has finite
$R$-dimension, the set of homogeneous elements of $A$ is an Ore set too and
the graded ring of fractions $Q^g$ is a gr-skewfield with $F^vQ$ inducing a
graded valuation in $Q^g$.
\end{prop}
\begin{proof}
Since $A_v$ is a domain and $G_{F^v}(A)=A_v[t,t^{-1}]$, $\sigma_v(S_0)$ is a
multiplicatively closed set, where $\sigma_v$ is the principal symbol map for
$F^{v}A$, and in fact $\sigma_v$ is a multiplicative map.  It follows that~:
$$
G_{F^v}(S^{-1}_0A)=\sigma_0(S_0)^{-1}G_{F^v}(A)=O_{\rm cl}(A_v)[t,t^{-1}]$$
the latter equality holds because the gradation is strong, hence localization
happens completely in degree zero.  Since the latter ring is a domain we may
apply Observation \ref{o2.1}. and following to $S^{-1}_0A$.  Note that in the finite
case we have $e=1$ because $\pi^n\Lambda\cap A_m=\pi^n\Lambda_m\neq
\pi^{n'}\Lambda_m$ for $n\neq n'$ in view of the Nakayama lemma.
\end{proof}
\sp
In general for given $A,FA$ (or $\wtA$) the construction of $\Lambda$,
$F\Lambda$ such that $\Lambda$ is $FA$-finite is not so easy, this problem is
related to the existence of discrete valuations having certain unramifiedness
properties.  In case the algebra is given by a finite number of generators and
finitely many relations between these, properties of so-called good reduction
will allow certain constructions of suitable $O_v$-orders.
\section{Positively Graded Connected Algebras}
A connected positively graded $K$-algebra is given as $A=K \oplus A_1\oplus
A_2 \oplus\ldots$, where $A_1$ is a finitely dimensional $K$-vector space $A$
is generated as a $K$-algebra by $A_1$.  We may view $A$ as a $K$-algebra
given by generators and homogeneous relations as follows~:
$$
0\longrightarrow{\cal R}\longrightarrow K\langle\Xul\rangle\mapright^{\pi} A\longrightarrow 0$$
where $K\langle\Xul\rangle$ is the free $K$-algebra on
$\Xul=\{X_1,\ldots,X_n\}$ and $\pi$ is defined by $\pi(X_i)=a_i$ where
$\{a_1,\ldots,a_n\}$ is a preselected $K$-basis of $A_1$.  The ideal ${\cal
R}$ is the ideal of relations.  By restricting $\pi$ to
$O_v\langle\Xul\rangle$ we obtain a graded subring $\Lambda$ of $A$ with
$\Lambda_0=O_v$ as follows~:
$$
0\mapright {\cal R}\cap O_v\langle\Xul\rangle \mapright O_v\langle\Xul\rangle
\mapright_{{\rm res}\pi}\Lambda\mapright 0$$
It is clear that $\pi$ maps $m_v\langle\Xul\rangle$ to $m_v\Lambda$ which is a
graded ideal of $\Lambda$.  So we have the following commutative diagram :
$$
\begin{diagram}
0 \rto & {\cal R}_v\rto & k_v\langle\Xul\rangle \rto & A_v\rto & 0\\
0 \rto &
{\cal R}\cap O_v\langle\Xul\rangle \rto\uto
\ar@{^{(}->}[d]
&
O_v\langle\Xul\rangle \rto_{{\rm res}\pi}\uto
\ar@{^{(}->}[d]
&
\Lambda \rto\uto
\ar@{^{(}->}[d]
& 0\\
0 \rto & {\cal R}\rto & K\langle\Xul\rangle \rto^{\pi} & A\rto & 0\\
\end{diagram}
$$
where $A_v=\Lambda/m_v\Lambda$ as usual and ${\cal R}_v={\cal R}\cap
O_v\langle\Xul\rangle\rangle+m_r\langle\Xul\rangle/m_r\langle\Xul\rangle$.
When ${\cal R}$ is generated by $\{p_1(\Xul),\ldots,p_d(\Xul)\}$ as a
two-sided ideal then, without loss of generality, we may assume that
$p_i(\Xul)\in O_v\langle\Xul\rangle$ (up to multiplying by a suitable
constant) but such that not all of them are in $m_v\langle\Xul\rangle$.
However the foregoing does not imply that ${\cal R}\cap
O_v\langle\Xul\rangle=O_v\langle\Xul\rangle p_1(\Xul)+\ldots
+O_v\langle\Xul\rangle p_d(\Xul)$, nor that ${\cal
R}_v=k_v\langle\Xul\rangle\pol_1(\Xul)+\ldots+
k_v\langle\Xul\rangle\pol_d(\Xul)$, where $\pol_i(\Xul)$ is the image of 
$p_i(\Xul)$ under reduction.
\begin{de}
With conventions and notation as before, we say that {\bf ${\cal R}$ (or $A$)
reduces well} or that $\Lambda$ defines a {\bf good reduction of $A$} whenever
${\cal R}_v$ is generated by the residues $\pol_i(\Xul),i=1,\ldots,d$, i.e.
whenever ${\cal R}\cap O_v\langle\Xul\rangle$ is generated by the $p_i(\Xul)$,
$i=1,\ldots,d$, as a two-sided ideal of $O_v\langle\Xul\rangle$.  Since $\pi$
is a graded morphism, the fact that ${\rm dim}_KA_1=n$ entails that ${\cal
R_1}=0$; it follows that ${\rm dim}_K(\Lambda/m_v\Lambda)_1=n$ but ${\rm
dim}_KA_n$ and ${\dim}_KA_{v,n}$ may be different from $n>1$.
\end{de}
\begin{prop}\label{p3.2}
Let $A=K[A_1]$ be a connected affine prime finite graded $K$-algebra and
$\pi:K\langle\Xul\rangle\r A$ a presentation of the $K$-algebra $A$ as in the
above diagram.  Then $F^vA$ is separated and $\cap$ is $FA$-finite if and only
if for all $n\in\Bbb{N}, {\rm dim}_KA_n={\rm \dim}_{k_v}A_{v.n}$.  Moreover, if
$A_v$ is a Goldie domain then $A$ is a domain and the $m_v\Lambda$-adic
filtration $F^vA$ is induced by a valuation filtration on the skewfield of
microfractions $\Delta$ of $A$.
\end{prop}
\begin{proof}
From Proposition \ref{p2.9} we retain that $\Lambda$ is $FA$-finite.  Suppose that
$A_v$ is a domain then we claim that $\Lambda$ and $A$ are domains too (we
cannot use Proposition \ref{p2.12} here because here $A$ is not necessarily an order
in a semisimple Artinian ring, in other words the Goldie ring property does
not follow from our assumptions unless we start from a Noetherian A) It will
be sufficient to check that there are no homogeneous zero-divisors.  Take
$a\in\Lambda_n, b\in \Lambda_m$ such that $ab=0$, say~: $a=\sum a_ix_i^{(n)},
b=\sum b_jx_j^{(m)}$, where $\{x^{(n)}_1,\ldots,x^{(n)}_{d_n}\}$ is an
$O_v$-basis for $\Lambda_n$. There exists $\mu_1,\mu_2\in K$ such that
$\mu_1a\in \Lambda_n,\mu_2b\in \Lambda_m$ but not all $\mu_1a_i$ in
$m_v\Lambda$, since is completely prime in $\Lambda$, either $\mu_1a\in m_v\Lambda$
or $\mu_2b\in m_v\Lambda$, a contradiction.  Now $F^vA$ is a strong filtration
and $G_{F^v}(A)=A_v[t,t^{-1}]$. Since $A_v$ is a Goldie domain, also
$A_v[t,t^{-1}]$ is a Goldie domain and it has a skewfield of fractions
$\Delta_v$ as well as a gr-skewfield of homogeneous fractions $\Delta^g_v$.
The multiplicative set $A-\{0\}$ has $\sigma_{F^v}(A-\{0\})=G_{F^v}(A)-\{0\}$
which is an Ore set because $A_v[t,t^{-1}]$ is a Goldie domain.  Put $S_0=\Lambda-\{0\}$, then $\sigma(S_0)=G_{F^v}(\Lambda)_--\{0\}$ where
$G_{F^v}(\Lambda)_-=\oplus_{n\le 0}G_{F^v}(\Lambda)_n$.  Clearly $\sigma(S_0)$
is an Ore set in $G_{F^v}(\Lambda)_-$ and also $\sigma(S_0)^{-1}G_{F^v}(\Lambda)_-=Q^g_{\rm cl}(G_{F^v}(A))_0$. Given $s\in S_0, a\in A$, the left Ore
condition for $\sigma(S_0$ yields $s'\in S_0,a'\in A$ such that $s'a-a's\in
m_v\Lambda$, say $s'a-a's=\pi^mb$ with $b\in \Lambda$ and $m\in\ $\d{$\Bbb{N}$}.
There are $s''\in S_0$ and $a''\in\Lambda$ such that $s''b-a''s\in
m_v\Lambda$, thus $(s''s')a-(s''a')s=\pi^m(a''s+\pi^{m'}y)$ for some
$y\in\Lambda$, $m'\neq 0$ in $\Bbb{N}$. Hence
$sa-(s''a'+\pi^ma'')s=\pi^{m+m'}$ and it follows that for every $p\ge q$, $S$
maps to an Ore set $S^{(p)}$ of $\Lambda/m^p_v\Lambda$.  Consequently the
microlocalization of $\Lambda$ at $S_0$ may be identified by 
$$\mathop{\rm lim}\limits_{\mapleft\limits_p}
(S^{(p)})^{-1}(\Lambda/m^p_r\Lambda)=Q^{\mu}(\Lambda)=Q^{\mu}({\cal
A}),$$
see \cite{AVVO} for more detail on microlocalization. Clearly $Q^{\mu}(A)$
is a skewfield and it has a strong filtration with associated graded ring
$\Delta^g_v$ that a domain and a graded skewfield. Applying Proposition 2.12
we may conclude that the filtration on $Q^{\mu}(A)$ is a discrete valuation
filtration.
\end{proof}
\sp
When considering a filtered $K$-algebra $A$ with a finite filtration $FA$, we
observe that there is an $n_0\in\Bbb{Z}$ such that for $n\le n_0$,
$F_nA=F_{n_0}A$.  Since we restricted attention to separated filtrations this
means that $F_nA=0$ for all $n\le n_0$ i.e. the filtration is left limited,
$F_{-1}A$ is a nilpotent ideal of $F_0A$.  Therefore, when dealing with finite
filtrations, it is not really restrictive to restrict attention to positively
filtered rings as we will do.
Moreover when domains have to be considered, $F_0A$ will be an algebraic field
extension of $K$ and so $O_v$ may be replaced by a discrete valuation ring of
$F_0A$ lying over $O_v\subset K$.  In other words we are lead to consider the
case of a positively filtered domain $K=F_0A\subset\ldots\subset
F_nA\subset\ldots\subset A$ and a discrete valuation ring $O_v$ of $K$ with an
$O_v$-order $\Lambda$ in $A$ such that $\Lambda\cap K=O_v$, $K\Lambda=A$,
equipped with the induced filtration $F\Lambda$ and $G_F(\Lambda)\subset
G_F(A)$ a graded $O_v$-order.  In the ``positive'' situation we have the
following lifting result.
\begin{prop}\label{p3.3}
If $G_F(\Lambda)$ is $G_F(A)$-finite then $\Lambda$ is $FA$-finite.
\end{prop}

\begin{proof} One easily establishes that $rk(F_q\Lambda)={\rm dim}_KF_qA$
by induction on $q$.  The case $q=0$ is trivial enough.  Assume that the
equality holds for $q-1$.  From $F_q\Lambda/F_{q-1}\Lambda=G(\Lambda)_q$ it
follows that~:
\begin{eqnarray*}
rk(F_q\Lambda)&=&
rk(F_{q-1}\Lambda)+{\rm \dim}_K(G(\Lambda)_q)={\rm \dim}_K(F_{q-1}A)+{\rm
\dim}_K(G(A)_q)\\
&=& {\rm \dim}_K(F_qA)\hskip8cm
\end{eqnarray*}
\end{proof}
\begin{cor}
Under the hypothesis of the foregoing Proposition 3.3. the results of Theorem
\ref{t2.6} are valid; in particular $G_{f^v}G_F(A)=G(A)_v[t,t^{-1}]=G_{f}G_{F^v}(A)$ where $f_nG_{F^v}(A)=F_nA_v[t,t^{-1}]$ is the filtration induced by $F$ in
$G_{F^v}(A)$ (this is a version of a general compatibility result for
arbitrary filtrations, cf. (\cite{MVO}, Proposition \ref{p2.4})).
\end{cor}
\begin{prop}\label{p3.5}
If $G(A)_v$ is a domain then also $G_F(A)$, $G_{F^v}(A)$ and $A$ are domains.
\end{prop}
\begin{proof}
Easy from the compatibility result for filtrations applied to $F^v$ and $F$,
i.e. $G_{f^v}G_F(A)=G_{f}G_{F^v}(A)$.
\end{proof}
\begin{cor}
If $G_F(\Lambda)$ is $G_F(A)$-finite then~:
\begin{enumrom}
\item
${\rm \dim}_{k_v}(A_{v,n})=rk G_F(\Lambda)_n={\rm \dim}_KG_F(A)_n$
\item
$\sum^n_{m=1}{\rm \dim}_{k_v}(A_{v,m})=\sum^n_{m=1}rk
G_F(\Lambda)_mrk(F_n\Lambda)\\\quad\\
={\rm dim}_KF_nA=\sum^n_{m=1}{\rm
\dim}_K(G(A)_m)$.
\end{enumrom}
\end{cor}

Assuming that $A=K[F_1A]$ then $A$ may be obtained as an epimorphic image of
the free $K$-algebra $K\langle\Xul\rangle$ in $\dim_KF_1A$-letters, say
$X_1,\ldots,X_d$, letting $\left\{x_1,\ldots,x_d\right\}$ be a $K$-basis for $F_1A$.
$$
\pi:K\langle X_1,\ldots,X_d\rangle\longrightarrow A,\quad X_i\longmapsto x_i,\quad i=1,\ldots,d$$
The filtration on $K\langle X_1,\ldots,X_d\rangle$ is the degree filtration
and this makes $\pi$ a strict filtered morphism in the sense of \cite{HVO1}.
Writing ${\cal R}={\rm Ker}\pi$, we have a {\bf strict exact} sequence of
filtered objects~:
$$
(*)\quad	0\longrightarrow {\cal R}\longrightarrow K\langle X_1,\ldots,X_d\rangle\longrightarrow A\longrightarrow0
$$
Strict exactness of (*) entails that by passing to Rees objects one obtains an
exact sequence of graded $K\langle X_1,\ldots,X_d\rangle$-modules:
$$
0\longrightarrow\wt{\cal R}\longrightarrow K\langle X_1,\ldots,X_d\rangle^{\sim}\longrightarrow \wtA\longrightarrow 0
$$
Again from strict exactness it follows that $G_{\cal F}(A)=G_F(A)$ where
$G_{\cal F}(A)$ is the associated graded of a $A$ as a filtered ${\cal
F}$-module, writing ${\cal F}=K\langle X_1,\ldots,X_d\rangle$.  From (*) we
thus derive an exact sequence in $G({\cal F})$-gr~:
$$
0\longrightarrow G({\cal R})\longrightarrow G({\cal F}) \longrightarrow G(A)\longrightarrow 0
$$
The filtration on ${\cal F}$ is exactly the gradation filtration it follows
that $G({\cal F})\cong {\cal F}$ and under this isomorphism $G({\cal R})$
corresponds to the ideal $\dot{\cal R}$ in ${\cal F}$ being the graded ideal
generated by the homogeneous components of highest degree in the homogeneous
decompositions of elements of ${\cal R}$.  The following is a version of
Theorem 2.13 in \cite{MVO}.
\begin{thm}\label{t3.7}
Assume that $A$ is given by finitely many generators and relations, where $FA$
is a described before Proposition \ref{p3.3}. If $G_F(A)$, which is a connected
positively graded $R$-algebra, reduces well with respect to $O_r$ then
$G_F(A)={\cal F}={\cal F}/\dot{\cal R}$, where ${\cal R}$ is generated as a
two-sided ideal by $p_1(\Xul),\ldots,p_s(\Xul)$ having
as homogeneous parts having highest degree $q_1(\Xul),\ldots,q_s(\Xul)$ that
generate $\dot{\cal R}$ as a two-sided ideal. Moreover $A$ reduces well at
$O_v$, in other words~:
$$
{\cal R}\cap O_v\langle \Xul\rangle=\sum_iO_v\langle
\Xul\rangle p_i(\Xul)O_v\langle\Xul\rangle$$
and $A_v$ is defined by the relations $p^v_1(\Xul),\ldots,p^v_s(X)$.\hfill\qed
\end{thm}
\sp
Foregoing theorem completes information about lifting properties of $G_F(A)$
to $A$ connected to the existence of valuation rings extending $O_v$ in either
$Q_{\rm cl}(A)$ if this exists (Noetherian or Goldie ring situation) or else
in a corresponding micro-localization $Q^{\mu}_{\rm cl}(A)$.  The finiteness
properties with respect to $\Lambda$ then provide the unramifiedness of the
extension of the valuation.  The latter unramified situation has been observed
in several independent interesting examples e.g.~:
\begin{enumrom}
\item
$\Delta(g)=Q_{\rm cl}(\mathcal{U}(g))$ for finite dimensional Lie algebras, Weyl
algebras $A_n({\Bbb{C}})$, cf. \cite{VOW}.
\item
Sklyanin algebras, cf. \cite{MVO}.
\item
Generalized gauge algebras including Witten algebras, cf. \cite{HVO1}.  The problem of
finding an extending noncommutative valuation has been reducted to finding an
$O_v$-order in an associated graded algebra having the finiteness property we
discussed and having a domain for its reduction.
\end{enumrom}
\section{Another Example: Generalization Weyl Algebras}
A generalized crossed product $A$ is a $\Bbb{Z}$-graded ring such that 
$A_i=A_0v_i$ is a free left $A_o$-module of rank one, and $v_0=l_A$ identifying 
$A_0$ as the subring $A_01_A$ in $A$.  Multiplication of $A$ is defined by: 
$$
av_ibv_j = a\sigma^i(b )c(i,j)v_{i+j}\ {\rm for}\  i,j \in\Bbb{Z},a\ {\rm and} 
\ b \in A_0$$ 
where $\sigma$ is an automorphism of $A_0$ and $c:\Bbb{Z}\times\Bbb{Z}\r\Bbb{Z} (A_0)$ is a 2-cocycle satisfying~: $c(i,j)c(i+j,k)=\sigma^i(c(j,k))c(i,j+k)$ 
for $i,j,k \in\Bbb{Z}$. A generalized Weyl algebra in the sense of 
(\cite{B1},\cite{B2},\cite{BVO}),\cite{BVO2}), is as before but now letting $A$ be generated over $A_0$ by two 
indeterminates $X=v_1$ and $Y=v_{-1}$ such that~: 
$$
\begin{array}{ll}
Xa = \sigma(a)X,\quad Ya = \sigma^{-1}(a)Y \quad {\rm for}\  a \in A_0\quad YX = a,\quad XY = \sigma(a)
\end{array}
$$
If $A_0=D$ is a commutative ring, e.g. a Dedekind domain, then these rings have 
now been extensively studied. Even over a Dedekind domain the class of 
generalized Weyl algebras contains many popular algebras: the first Weyl 
algebra and its quantum deformation, the quantum plane, the quantum 
2-dimensional sphere, $\mathcal{U}(sl_2)$ and its quantum version $\mathcal{U}_q(sl_2)$ Witten's 
first deformation and Woronowicz's deformation, the quantum Heisemberg algebra, 
the Virasaro algebra.  We write $D(a,\sigma)$ for generalized Weyl algebra as 
above with $A_0=D$. 
\sp
We now consider $\K\subset D$ a fixed base field invariant under $\sigma$. 
Write $D^{\langle\sigma\rangle}$ for the invariant algebra with respect to 
$\sigma$.  We may restrict attention to affine $\Bbb{K}$-algebra $D$ but the 
results can be generalized to the consideration of Noetherian integrally closed domains (localization at height one ideals then yields discrete valuation 
rings).  Localizing $D(\sigma,a)$ at $D-\{0\}$ yields 
$\Bbb{K}\otimes_DD(\sigma,a)\simeq\Bbb{K}[t,t^{-1}\sigma]$, where $\Bbb{K}=
Q_{cl}(D)$, $\sigma$ the extended automorphisms of $\Bbb{K}$. 
\begin{lem}
If $P$ is a $\sigma$-invariant prime ideal of $D$ then $D(\sigma,a)P$ is a 
two-sided ideal such that $D((\sigma,a)/D(\sigma,a)P$ is of type $\Dol
(\ol{\sigma},\aol)\Pol$ where $\Dol=D/P,\sigma$ is induced by $\sigma$ on 
$\Dol$ and $\aol=a$ mod $P$ ($\aol=0$ is allowed).
\end{lem} 
\begin{proof} 
The $\sigma$-invariance of $P$ yields that $D(\sigma,a)P$ is two-sided. Since 
$D$ is Dedekind, $\Dol$ is a field.  If $a\not\in P$ then 
$\Dol(\ol{\sigma},\aol)$ is again a generalized Weyl algebra and a domain. 
If $a \in P$, then $\Dol(\ol{\sigma},a)$ is not a domain. 
If $D=O_v\subset\Bbb{K}$ the maximal ideal is necessarily $\sigma$-invariant 
and, $\Dol(\ol{\sigma},\aol)=\ol{\Bbb{K}}(\ol{\sigma},\aol)$.  
If $\aol\neq 0$ is necessarily a unit of $\ol{\Bbb{K}}$. 
\end{proof}
\sp
More generally, If $P$ is $\sigma$-invariant in $D$ then $D-P$ is also 
$\sigma$-invariant hence an Ore set in $D(\sigma,a)$. Localizing $D(\sigma,a)$ 
at $D-P$ then yields $D(\sigma,a)_P=D_P(\sigma,a)$.  If $P\neq 0$ then $D_P$ is 
a discrete valuation ring of $\Bbb{K}$ and $D_P(\sigma,a)$ is a gr-valuation 
ring in $\Bbb{K}[t,t^{-1},\sigma]$ (the latter being a graded-skewfield).
The corresponding valuation filtration on $\Bbb{K}[t,t^{-1},\sigma]$ is compatible with
$\Bbb{Z}$-grading and the associated graded ring for the valuation filtration 
is exactly $\ol{\Bbb{K}}[t,t^{-1},\sigma]$.  It follows that $D_P(\sigma,a)$ is 
an intersection of $\Bbb{K}[t,t^{-1},\sigma]$ and a discrete valuation on 
$Q_{cl}(D(\sigma,a))=\Bbb{K}[t,\sigma]$.  So we have proved. 

\begin{prop} 
If $P$ is a $\sigma$-invariant prime ideal of $D$ such that $a \not\in P$ then 
$P$ determines a noncommutative valuation of $Q_{cl}(D(\sigma,a))\cong \Bbb{K}
(t,\sigma)$ with ring $\Lambda_P$ say, and maximal ideal $w$, such that 
$D(\sigma,a)=\Lambda_P\cap\Bbb{K}[t,t^{-1},\sigma]$, and $P=w\cap D(\sigma,a)$. 
The residue skew field of this discrete valuation is 
$Q_{cl}(\Dol(\ol{\sigma},\aol))=\ol{\Bbb{K}}[t,\ol{\sigma}]$. 
\end{prop}
\sp
Note that for the Weyl algebra 
$$A_1(\Bbb{C})=\Bbb{C}\left[X,Y\right]/\langle XY-YX-1\rangle,\quad D_1=\Bbb{C}[X,Y],\quad \sigma(XY)=XY+1$$ 
there are not nontrivial 
$\sigma$-invariant prime ideals in $D_1$.  In a sense the prime at $\infty$ is an invariant prime (corresponding to $\Bbb{C}[(XY)^{-1}]_{(xy)^{-1}})$ 
and it is the valuation ring in $D_1(\Bbb{C})=Q_{cl}(A_1(\Bbb{C}))$ 
corresponding to the quotient filtration of the Bernstein filtration on 
$A_1(\Bbb{C})$ that represents this prime at $\infty$ (we refer to \cite{VOW} for 
some  results on valuations of $D_1(\Bbb{C}))$.  
\sp
Look at $D_1$, the coordinate ring of curve $C$ in affine $n$-space over 
$\Bbb{K}$, in particular $\Bbb{K}$ is algebraically closed in the field of 
fractions of $D_1,{\bf K}$ say. 
\sp Let $O_v\subset\Bbb{K}$ be such that $C$ has good reduction at $O_v$ i.e. 
the reduced equations of $C$ define a nonsingular curve over the residue field 
$\ol{\bf K}$, or equivalently $\ol{D}= D_{O_v}/m_vD_{O_v}$ is a Dedekind domain,
where $D_{O_v}=O_v[X_1,\ldots,X_n]/I\cap O_v[X_1,\ldots,X_n],I$ 
is an ideal of $C$.  If $a\in D_{O_v}-m_vD_{O_v}$ then $D_{O_v}(\sigma,a)= 
\Lambda\subset D(\sigma,\Lambda)$.  It is clear that $m_vD_{O_v}$ is 
$\sigma$-invariant, therefore the left ideal $D_{O_v}(\sigma,a)m_v$ is 
two-sided and we have~: 
\begin{lem} 
 as before : $D_{O_v}(\sigma,a)/m_vD_{O_v}(\sigma,a)\cong 
\Dol(\ol{\sigma},\aol)$. 
If the curve given by $D$ over $\Bbb{K}$ has good reduction at $O_v$ then for 
$a \not\in m_vD$ we obtain an $m_v$-adic filtration on $D(\sigma,a)$ with 
$F_nD(\sigma,a)=D_{O_v}(\sigma,a)\pi^{-n}$ where $m_v=(\pi)$ for $n\in\Bbb{Z}$ 
such that the associated graded ring is $\Dol(\ol{\sigma},\aol)[t,t^{-1}]$,
hence a domain. 
\end{lem}
\begin{rem}
The restriction to $O_v$ defining good reduction for $C$ can be avoided. 
In the above result it is only important to have $\Dol(\ol{\sigma},\aol)$ to be 
a domain and the assumption on a shows that it is enough to have 
$\Dol(\ol{\sigma},\aol)$ to be a domain.  However from the point of view of 
the theory of generalized Weyl algebras it is nice to have a residue algebra 
again being a generalized Weyl algebra of the same type.  Therefore ``good 
reduction'' is an interesting condition.  Combining this with the results of
section 1, we may phrase all this as follows. 
\end{rem}
\begin{thm} 
Let $\Bbb{K}\subset D^{\langle\sigma\rangle}\subset D$ where $D$ is the 
coordinate ring of a nonsingular curve $C$ in $n$-space.  Let $m_b\subset
O_v\subset \Bbb{K}$ define a discrete valuation of $\Bbb{K}$ such 
that $a\in D_{O_v}-m_vD_{O_v}$.  The filtration $FD(\sigma,a)$ defined by $F_nD(\sigma,a)$ extends to a valuation filtration on the skewfield $Q_{cl}(D(\sigma,a))\cong\Bbb{K}(t,\sigma)$ with residue skewfield 
$A_{cl}(\Dol(\ol{\sigma},\aol))$. 
\end{thm}
If $C$ has good reduction at $O_v$ then $\Dol(\ol{\sigma},\aol)$ is a 
generalized Weyl algebra over the Dedekind domain $\Dol$. 
If $\sigma$ has infinite order then $D^{\langle\sigma\rangle}=\Bbb{K}$. 
\begin{proof}
Only the final statement has not yet been fully established.  If $D^{\langle
\sigma\rangle}$ is not algebraic over $\Bbb{K}$ then ${\bf K}$ must be 
algebraic over $Q_{cl}(D^{\langle\sigma\rangle})={\bf
K}^{\langle\sigma\rangle}$.  Since $D$ is affine over $\Bbb{K},{\bf K}$ is 
finitely generated as a field over $\Bbb{K}$ hence over ${\bf
K}^{\langle\sigma\rangle}$.  It follows that $[{\bf K}:{\bf K}]\le\infty$  but 
then $\langle\sigma\rangle$ is a finite group, a contradiction. 
\end{proof}
\begin{cor}
Certain discrete valuations of the base field $\Bbb{K}$ extend to 
noncommutative discrete valuations (unramified extension) on the skewfield 
of fractions of quantum enveloping algebras, the quantum plane, the quantum 
$O_{q^2}$ of $so(\Bbb{K},3)$ \cite{Sm}, the quantum Heisemberg algebra \cite{Ma}, generalized 
gauge algebra of \cite{LVO}.  The condition on discrete valuation of the base field 
is given in terms of good reduction of some constant e.g. $q$ or $a$.  For 
example in case of the quantum Weyl algebra $A=\Bbb{K}[t](\sigma,t)$ where 
$\sigma(t)=q^{-1}(t-1)$ it is clear that $\ol{\sigma}$ explodes when one allows 
an $O_v$ of $\Bbb{K}$ containing $q$ in $m_v$.
\end{cor}

\end{document}